\documentclass[a4paper,11pt]{article}
\usepackage{citesort}
\newtheorem{theorem}{Theorem}[section]
\newtheorem{lemma}[theorem]{Lemma}

\newtheorem{prop}[theorem]{Proposition}
\newtheorem{defi}[theorem]{Definition}
\newtheorem{coro}[theorem]{Corollary}

\newcommand{\eproof}{{\vrule height8pt width5pt depth0pt}\vspace{3mm}}

\newcommand{\eg}{\varepsilon}
\newcommand{\llg}{\lambda}
\newcommand{\ag}{\alpha}
\newcommand{\bg}{\beta}

\newcommand{\Og}{\Omega}
\newcommand{\pdh}{\partial}
\newcommand{\RR}{{\rm I\kern -1.6pt{\rm R}}}
\newcommand{\mO}{\mathcal{O}}
\newcommand{\mD}{\mathcal{D}}

\newcommand{\mF}{\mathcal{F}}
\newcommand{\mP}{\mathcal{P}}

\newcommand{\mR}{\mathcal{R}}

\newcommand{\mT}{\mathcal{T}}

\title{On singularities and instability of reconstruction in thermoacoustic tomography}
\author{Linh Nguyen}

\begin{document} \maketitle

\begin{abstract} We consider the problem of thermoacoustic tomography (TAT), in which one needs to reconstruct the initial value of a solution of the wave equation from its value on an observation surface. We show that if some geometric rays for the equation do not intersect the observation surface, the reconstruction in TAT is not H\"{o}lder stable. \end{abstract}


\section{Introduction}\label{intro}Thermoacoustic tomography (TAT) is a hybrid imaging method that combines the significant contrast of electromagnetic and high resolution of ultrasound imaging. A biological object is irradiated by a brief electromagnetic (EM) pulse in visible light or radiofrequency range. Some fraction of the EM energy is absorbed by the tissues. Since EM absorption is much higher in tumors, knowing the distribution $a(x)$ of the absorbed energy would provide some valuable diagnostic information. The absorption causes thermoelastic expansion in the tissues, and thus triggers a pressure (ultrasound) wave $p(x,t)$ propagating through the body. The pressure is then measured on an observation surface $S$ completely or partially surrounding the object. The initial pressure $f(x) = p(x,0)$ is roughly proportional to the EM energy absorption $a(x)$. One thus concentrates on the recovery of $f(x)$ instead of $a(x)$.

The commonly accepted model of TAT is (e.g.,\cite{DSK,Tam}) \begin{eqnarray} \label{E:wave} \left \{\begin{array}{l} p_{tt}(x,t) = c^2(x) \bigtriangleup  p(x,t),\quad  x \in \RR^n,\quad t \geq 0, \\ p(x,0) =f(x), \quad p_t(x,0) = 0,\\p(y,t) =g(y,t), y \in S, t \geq 0. \end{array} \right.\end{eqnarray} Here $c(x)$ is the the ultrasound speed at location $x$, $g(y,t)$ is the pressure measured at location $y \in S$ and time $t$. We assume that the sound speed $c(x)$ is smooth and there are $C_0,c_0>0$ such that $c_0 \leq c(x) \leq C_0$ for all $x \in \RR^n$. The main problem of TAT is: given $c(x)$, reconstruct $f$ from $g$. Significant progresses have been made in mathematics of TAT. The readers are referred to \cite{AKKun,XW06,KKun,FPR} for thorough discussion and relevant references.

In this text, we concentrate on the stability analysis of the reconstruction of $f$ in a region of interest $\Omega$, which is a bounded open domain in $\RR^n$, from the measured data $g$ on the observation surface $S$. This issue has been partially addressed in a number of papers \cite{XWAK04,XWAK08,KunOpen,FPR,HKN,US,HTRE,Pal,LQ00}.

We first recall the analysis for the case of constant speed $c(x)=c$. It was shown that if the so-called \emph{visibility} condition is satisfied, the reconstruction is only mildly unstable, similarly to the inversion of the standard Radon transform (see, e.g. \cite{KunOpen}, for more discussion and references). The visibility condition requires that for each point $x$ in $\Omega$, every straight line passing through $x$ intersects the observation surface $S$. For instance, if $S$ is a closed hypersurface surrounding $\Omega$ then the visibility condition is satisfied.

Let us now consider the case of variable speed $c(x)$.  If $\Og$ is enclosed by $S$, it was argued and demonstrated in numerics in \cite{HKN} that if there are some geometric rays trapped inside $\Og$, then the singularities of $f$ that lead to these rays cannot be stably reconstructed. On the other hand, the authors of \cite{US} recently proved that under visibility condition, which can be defined for variable speed using geometric rays instead of straight lines, the reconstruction is Lipschitz stable. In this paper, we prove a complementary result, which shows that the reconstruction is not H\"{o}lder stable if the visibility condition does not hold. The visibility condition can be violated when either the data is incomplete ($S$ does not completely surround $\Og$, see, e.g., \cite{KunOpen}), or the trapping phenomenon occurs (e.g., \cite{HKN}). These two cases have the same instability, and we do not distinguish between them. Although we do not use this in the text, it should be mentioned that the Lipschitz instability can be obtained from the general framework proposed in \cite{USL}.

The text is organized as follows. In section \ref{S:SW}, we recall the notion of wavefront set, whose propagation is the central issue of stability analysis in TAT (as well as other types of tomography, e.g., \cite{LQ00,QT93,FLU,GrUDuke,GrUCon}). We then present our instability result in section \ref{S:inst}.

\section{Singularities and Wavefront set}\label{S:SW}
In tomography, singularities are usually related to sharp details, for instance the boundaries of objects, jumps in densities, or interfaces between tissues. In many cases, singularities (rather than the exact image) are of the main interest (e.g., they are the objects of reconstruction in local tomography, e.g. \cite{KLM}). We will see that propagation of singularities also plays an important role in the stability analysis of reconstruction in TAT.

We denote by $\mD(\RR^n)$ and $\mD'(\RR^n)$ the standard spaces of test functions and distributions on $\RR^n$. We now recall the definition of wavefront set (e.g., \cite{STGuide}):
\begin{defi} Let $f \in \mD'(\RR^n)$, $x_0 \in \RR^n$, and $\xi_0 \in \RR^n \setminus \left\{0\right\}$. Then $f$ is microlocally smooth at $(x_0,\xi_0)$ if there is a function $\varphi \in C_0^\infty(\RR^n)$ satisfying $\varphi(x_0) \neq 0$ and an open cone $V$ containing $\xi_0$, such that $\mF (\varphi f)$ is rapidly decreasing in $V$. That is, for any $N>0$  there exists a constant $C_N$ such that
$$|\mF(\varphi f)(\xi)| \leq C_N \left<\xi\right>^{-N}, \mbox{ for all } \xi \in V.$$ Here $\mF$ is the Fourier transform and $\left<\xi\right>= \left(1+|\xi|^2 \right)^{1/2}$. The wavefront set of $f$, denoted by $WF(f)$, is the complement of the set of all $(x_0,\xi_0) \in \RR^n \times (\RR^n \setminus 0)$ where $f$ is microlocally smooth.
\end{defi} For example, if $f$ is the characteristic function of an open set $\Og$ with smooth boundary $\pdh \Og$, then $(x_0,\xi_0) \in WF(f)$ if and only if $x_0 \in \partial \Og$ and $\xi_0$ is perpendicular to the tangent plane $T_{x_0}\partial \Omega$ of $\partial \Og$ at $x_0$.

Propagation of the wavefront set of the solution $p(x,t)$ for equation (\ref{E:wave}) can be precisely described (e.g., \cite{NiLec,STGuide,US,Stei}). Since the initial velocity $p_t(x,0)$ is zero, each element $(x_0,\xi_0) \in WF(f)$ propagates in two opposite directions $\xi_0$ and $-\xi_0$. Let us consider the bi-characteristics $(x_\pm(s),t_\pm(s),\xi_\pm(s),\tau_\pm(s)),$ which are the solutions of the following Hamiltonian systems:
\begin{eqnarray*}\left \{\begin{array}{l} \dot{x}(s) = - c^2(x) \xi,~\dot{t}(s)= \tau, \\ \dot{\xi}(s) =\frac{1}{2} \nabla c^2(x) |\xi|^2,~\dot{\tau}(s) =0 \\ (x(0),t(0),\xi(0),\tau(0))=(x_0,0,\mp \xi_0,c(x_0)|\xi_0|).\end{array} \right.\end{eqnarray*} Since $c(x)$ is smooth and $0< c_0 \leq c(x) \leq C_0$, these bi-characteristics are well defined on $s \in \overline{\RR}_+$, and $(x_\pm(s),t_\pm(s),\xi_\pm(s),\tau_\pm(s)) \in WF(p)$ for all $s$. We denote the (x,t)-projections of these bicharacteristics by $\mR_{+}(x_0,\xi_0)$ and $\mR_{-}(x_0,\xi_0)$. The following result can be easily obtained by basic tools of microlocal analysis (e.g., \cite{STGuide}):
\begin{theorem}\label{T:smooth} Consider equation (\ref{E:wave}). Let $\mO$ be an open subset of $\RR^n \times \overline{\RR}_+$. Assume that for all $(x,\xi) \in WF(f)$, $\mR_+(x,\xi)$ and $\mR_-(x,\xi)$ do not intersect  $\mO$. Then $p \in C^\infty(\mO)$.\end{theorem}
Let $\mR_x(x,\xi)$ be the x-projection of $\mR_+(x,\xi) \cup \mR_-(x,\xi)$. $\mR_x(x,\xi)$ is, indeed, a connected smooth curve in $\RR^n$, which we call a (geometric) ray. The following result is a simple consequence of Theorem \ref{T:smooth}:
\begin{coro}\label{C:smooth} Consider equation (\ref{E:wave}). Let $V$ be an open subset of $\RR^n$. Assume that for all $(x,\xi) \in WF(f)$, the rays $\mR_x(x,\xi)$ do not intersect  $V$. Then $p \in C^\infty(V \times \overline{\RR}_+)$.\end{coro}

\section{Instability of reconstruction in thermoacoustic tomography}\label{S:inst}
Let us return to the main equation of TAT: \begin{eqnarray} \label{E:WaveInst} \left \{ \begin{array}{l} p_{tt}(x,t) = c^2(x) \bigtriangleup p(x,t),\quad \forall x \in \RR^n,\quad t>0, \\ p(x,0) =f(x), \quad  p_t(x,0) = 0.\end{array} \right.\end{eqnarray} Let $S$ be a closed piece of a hypersurface in $\RR^n$, $T>0$, $\Gamma=S \times [0,T]$, and $g$ be the restriction of $p$ on $\Gamma$. We define the linear operator \begin{eqnarray*} \mT: L^2(\Og) &\longrightarrow & \mD'(\Gamma) \\ f &\longmapsto& g.\end{eqnarray*} Here, we identify $L^2(\Og)$ with the subspace of $L^2(\RR^n)$ containing functions supported in $\overline{\Og}$. We assume that $\mT$ is injective and consider the stability problem of the reconstruction of $f$ from $g$. Let $U$ and $V$ be open sets in $\RR^n$ such that $U \subset \Og$ and $S \subset V$.
\begin{theorem}\label{T:holder} Assume that there exists a nonzero vector $\xi_0 \in \RR^n \setminus 0$ such that for all $x \in U$, the rays $\mR_x(x,\xi_0)$ do not intersect $V$. Then the reconstruction of $f$ from $g$ is not H\"{o}lder stable. That is, there do not exist $\mu>0$, $\delta>0$, $s_0,s_1 \geq 0$, and $C>0$ such that $$\|f\|_{L^2(\Og)} \leq C \|\mT f\|^{\mu}_{H^{s_0}(\Gamma)}, \mbox{
for all } f \in H^{s_1}(\Og) \mbox{ satisfying }\|f\|_{H^{s_1}(\Og)} \leq \delta.$$ \end{theorem}

In order to prove this theorem, we need an auxiliary result. Without loss of generality, we can assume that $0 \in U$ and $\xi_0 = (0,..,0,1)$. Since $U$ is open, there is an open set $U_0 \subset \RR^{n-1}$ and $\eg>0$ such that $\overline{U}_0 \times \bar{I} \subset U$, where $I = (-\eg,\eg)$. Let us fix a nonzero function $f_0 \in C^\infty_0(U_0)$. For each $x \in \RR^n$, we write $x=(x',x_n)$, where $x' \in \RR^{n-1}$ and $x_n \in \RR$. We now consider $$X = \{f \in L^2(U): f(x) = f_0(x') h(x_n), h \in L^2(\RR), supp(h) \subset \overline{I}\}.$$ Then, $X$ is an infinite dimensional closed subspace of $L^2(\Og)$.
\begin{lemma}\label{L:auxi}  For all $s \geq 0$, $\mT$ induces a linear bounded operator $$\mT|_X: X \longrightarrow H^s(\Gamma).$$ \end{lemma}
In what follows, we first prove Lemma \ref{L:auxi}, and then show how to obtain Theorem \ref{T:holder} from it.

\subsection{Proof of Lemma \ref{L:auxi}}
Let $f \in X$, then $WF(f) \subset U \times \{\llg \xi_0: \llg \neq 0\}$. Since $\mR_x(x,\llg \xi_0)=\mR_x(x,\xi_0)$, one has $\mR_x(x,\xi)= \mR_x(x,\xi_0)$ for all $(x,\xi) \in WF(f)$. Assuming the condition in Theorem \ref{T:holder}, one deduces that $\mR_x(x,\xi)$ do not intersect $V$ for all $(x,\xi) \in WF(f)$.

We consider equation (\ref{E:WaveInst}) and let $\mP(f)=p|_{V\times \overline{\RR}_+}$. Due to Corollary \ref{C:smooth}, $\mP|_{X}$ is a linear operator from $X$ to $C^\infty(V \times\overline{\RR}_+)$. Let $V_0$ be an open set in $\RR^n$ such that $S \subset V_0 \subset \overline{V}_0 \subset V$. Then, for all $s \geq 0$ and $T>0$, $\mP$ induces a linear operator: $$\mP|_{X}:X \rightarrow H^s(V_0 \times [0,T]).$$

We now prove that this operator is bounded. Indeed, we first show that $\mP:L^2(\Og) \rightarrow L^2(V_0 \times [0,T])$ is bounded. Let $u(.,t) = \int_0^t p(.,s)ds$, then, from (\ref{E:WaveInst}), $u$ solves the equation
\begin{eqnarray}\label{waveP} \left \{\begin{array}{l} u_{tt}= c^2(x) \Delta u(x,t), ~ x \in \RR^n, ~ t>0, \\ u(x,0)=0, ~ u_t(x,0) = f(x).
\end{array}\right. \end{eqnarray}  Denoting $E(t)=\|c^{-1}u_t(.,t)\|^2_{L^2(\RR^n)} + \|\nabla u(.,t)\|^2_{L^2(\Og)},$  one obtains $E$ is independent of $t$, due to the conservation of energy (e.g., \cite{CH2}). That is, $E(t)=E(0)$ for all $t \in \RR_+$. Since $u_t(.,t) = p(.,t)$, $u(.,0)=0$, and $u_t(.,0)=f$, one derives $$\|c^{-1}p(.,t)\|^2_{L^2(\RR^n)} + \|\nabla u(.,t)\|^2_{\RR^n} = \|c^{-1} f\|_{L^2(\RR^n)}.$$ Using the inequalities $0< c_0 \leq c(x) \leq C_0$, one deduces that there is a constant $A>0$ satisfying $$\|p(.,t)\|_{L^2(\RR^n)} \leq A \|f\|_{L^2(\RR^n)}, \quad \mbox{for all } f \in L^2(\Og).$$ Hence, $$\|p\|_{L^2(V_0\times [0,T])} \leq \|p(.,t)\|_{L^2(\RR^n \times [0,T])} \leq AT \|f\|_{L^2(\RR^n)},$$ which proves the boundedness of $\mP:L^2(\Og) \rightarrow L^2(V_0 \times[0,T])$.

Since $H^s(V_0 \times [0,T])$ is continuously embedded into $L^2(V_0 \times [0,T])$, applying Propositions \ref{P:cont} (see appendix), one concludes that $\mP|_X: X \rightarrow H^s(V_0 \times [0,T])$ is bounded.

We are now ready to prove Lemma \ref{L:auxi}. Since the restriction $\mR(p) = p|_{\Gamma}$, as a linear operator $\mR: H^s(V_0 \times [0,T]) \longrightarrow H^{s-1/2}(\Gamma)$, is bounded for any $s > \frac{1}{2}$, we have $$\mT|_X= \mR \circ \mP|_X: X \longrightarrow H^s(\Gamma)$$ is bounded for any $s \geq 0$.

\subsection{Proof of Theorem \ref{T:holder}}

We first recall some facts concerning singular value decomposition (e.g., \cite{ETFunct,GKr}). Let $H_1,H_2$ be Hilbert spaces and $A:H_1 \rightarrow H_2$ be a bounded injective operator. Let $A^*$ be the adjoint operator of $A$ and $B=A^*A$. Then $B: H_1 \rightarrow H_2$ is a positive definite bounded operator. Let us denote by $\{s_j^2\}$ the eigenvalues of $B$ and by $\{e_j\}$ the corresponding unit norm eigenvectors. Then $\{e_j\}_j$ is an orthogonal basis of $H_1$. Denoting $f_j=\frac{1}{s_j} A(e_j)$, it is simple to show $\{f_j\}_j$ is an orthonormal set in $H_2$. If $range(A)$ is dense in $H_2$, then $\{f_j\}_j$ is an orthonormal basis of $H_2$.

For an operator $A$, we denote by $\{s_j(A)\}_j$ the set of above $s_j$, which are chosen to be positive and decreasing. Then $\{s_j(A)\}_j$ are called s-numbers (or singular values) of $A$. We will need the following asymptotic behavior of s-numbers:
\begin{lemma} \label{sing} (e.g., \cite[page 119]{ETFunct}) \begin{enumerate} \item Let $s_1>s_2$ and $i_1: H^{s_1}(\Gamma) \longrightarrow H^{s_2}(\Gamma)$ be the natural embedding. There exists a constant $c_1>0$ independent of $j$, such that $s_j(i_1) \leq c_1 j^{\frac{s_2-s_1}{n}}$. \item Let $s>0$, $\eg>0$, and $i_2: H^{s}_0(-\eg,\eg) \longrightarrow L^2(-\eg,\eg)$ be the natural embedding. There exists a constant $c_2>0$ independent of $j$, such that $s_j(i_2) \geq c_2 j^{-s}$.\end{enumerate} \end{lemma}

{\bf Proof of Theorem \ref{T:holder}} Suppose that there exist $\mu>0$, $\delta>0$, $s_0,s_1 \geq 0$, and $C>0$ such that \begin{eqnarray*}\|f\|_{L^2(\Og)} \leq C \|\mT f\|^{\mu}_{H^{s_0}(\Gamma)} \end{eqnarray*} for all $f \in H^{s_1}(\Og)$ satisfying $\|f\|_{H^{s_1}(\Og)} \leq \delta$.  Then for any $f \in H^{s_1}(\Og)$, applying this inequality to $\frac{\delta f}{\|f\|_{H^{s_1}(\Og)}}$, we have $$\frac{\|f\|_{L^2(\Og)}}{\|f\|_{H^{s_1}(\Og)}} \leq C_1 \left(\frac{\|\mT f\|_{H^{s_0}(\Gamma)}}{\|f\|_{H^{s_1}(\Og)}}\right)^{\mu},$$ where $C_1$ is a constant independent of $f$. That is, \begin{equation}\label{E:l2in} \|f\|_{L^2(\Og)} \leq C_1 \|\mT f\|_{H^{s_0}(\Gamma)}^{\mu}
\|f\|_{H^{s_1}(\Og)}^{1-\mu},\end{equation} for all $f \in H^{s_1}(\Og)$. Let us consider the following subspace of space $X$ from Lemma \ref{L:auxi}:
$$X_{s_1} = \{f : f(x) = f(x',x_n)=u_0(x') h(x_n), h \in H^{s_1}_0(I)\}.$$ We then have $X_{s_1} \subset X \cap H^{s_1}(\Og)$. We now prove that (\ref{E:l2in}) does not hold true for all $f \in X_{s_1}$. Indeed, due to Lemma \ref{L:auxi}, for any $f\in X_{s_1} \subset X$, one has \begin{equation}\label{tff}
\|\mT f\|_{H^s(\Gamma)}^\mu \leq C_s\|f\|_{L^2(\Og)}^\mu.\end{equation} Here, $s \geq 0$ and $C_s$ is a general constant depending on $s$. Combing (\ref{E:l2in}) and (\ref{tff}), we have $$\|Tf\|_{H^s(\Gamma)}^\mu
\|f\|_{L^2(\Og)}^{1-\mu} \leq C_s \|f\|_{L^2(\Og)} \leq C_s \|\mT f\|_{H^{s_0}(\Gamma)}^\mu\|f\|_{H^{s_1}(\Og)}^{1-\mu}.$$
That is, \begin{equation}\label{e1e} \left(\frac{\|f\|_{L^2(\Og)}}{\|f\|_{H^{s_1}(\Og)}}\right)^{1-\mu} \leq C_s \left(\frac{\|\mT f\|_{H^{s_0}(\Gamma)}}{\|\mT f\|_{H^s(\Gamma)}}\right)^\mu .\end{equation}
Since $f(x)= f_0(x')h(x_n)$, direct computations show \begin{eqnarray*}\|f\|_{L^2(\Og)} &=& \|f_0\|_{L^2(\RR^{n-1})} \|h\|_{L^2(I)},\\ \|f\|_{H^{s_1}(\Og)} &\leq& \|f_0\|_{H^{s_1}(\RR^{n-1})}\|h\|_{H^{s_1}(I)}.\end{eqnarray*}
Inequality (\ref{e1e}) gives \begin{equation}\label{htf}\left(\frac{\|h\|_{L^2(I)}}{\|h\|_{H^{s_1}(I)}}\right)^{1-\mu} \leq C_s \left(\frac{\|\mT f\|_{H^{s_0}(\Gamma)}}{\|\mT f\|_{H^s(\Gamma)}}\right)^\mu,
\mbox{ for all } h \in H^{s_0}_0(I).\end{equation} Let $s>s_0$, we now prove that this equality implies $[s_j(i_2)]^{1-\mu} \leq C_s [s_j(i_1)]^{\mu}$ for all $j$, where $i_1,i_2$ are the natural embeddings $i_1: H^{s}(\Gamma) \hookrightarrow H^{s_0}(\Gamma)$ and $i_2: H^{s_1}_0(I)\hookrightarrow L^2(I)$. Indeed, let $\{h_j\}_j$ be the orthonormal basis of $H^{s_1}_0(I)$ such that $$\|h_j\|_{L^2(I)}= s_j(i_2) \|h_j\|_{H^s(I)}.$$
Let $\ag= (\ag_1,...,\ag_j) \in \RR^j$ and $h= \sum_{k \leq j} \ag_k h_k$. Since $\{h_j\}_j$ are also orthogonal in $L^2(I)$,  \begin{eqnarray*}\|h\|^2_{L^2(I)}&=& \sum_{k \leq j} |\ag_k|^2 \|h_k\|^2_{L^2(I)} =  \sum_{k \leq j} |\ag_k|^2 s^2_k(i_2) \|h_k\|^2_{H^{s_1}(I)} \\ &\geq& s_j^2(i_2) \sum_{k \leq j} |\ag_k|^2  \|h_k\|^2_{H^{s_1}(I)} = s_j^2(i_2)\|h\|^2_{H^{s_1}(I)}. \end{eqnarray*}
That is, \begin{equation}\label{E:hin} \|h\|_{L^2(I)} \geq s_j(i_2)\|h\|_{H^{s_1}(I)}. \end{equation}

Now, consider $\{g_j\}_j$ to be the orthonormal basis of $H^{s}(\Gamma)$ such that $$\|g_j\|_{H^{s_0}(\Gamma)}=s_j(i_1)\|g_j\|_{H^{s}(\Gamma)}.$$
Fixing $j$, we can choose $\ag =(\ag_1,...,\ag_j) \neq 0$ such that, for $h=\sum_{k \leq j} \ag_k h_k$ and $f(x)=f_0(x') h(x_n)$, $\mT (f)$ belongs to the orthogonal complement of  $span\{g_1,...,g_{j-1}\}$. Since the embedding $i_1: H^s(\Gamma)\rightarrow H^{s_0}(\Gamma)$ has dense range, $\{g_k\}_k$  is also an orthogonal basis of $H^{s_0}(\Gamma)$. One then obtains $\mT f= \sum_{k \geq j} \bg_k g_k$ in both $H^{s_0}(\Gamma)$ and $H^s(\Gamma)$. Therefore, $$ \|\mT f\|^2_{H^{s_0}(\Gamma)}= \sum_{k \geq j} \bg_k^2 \|g_k\|^2_{H^{s_0}(\Gamma)}= \sum_{k \geq j} \bg_k^2 s^2_k(i_1) \|g_k\|^2_{H^{s}(\Gamma)} \leq s^2_j(i_1) \|\mT f\|^2_{H^{s}(\Gamma)}.$$
That is, \begin{equation} \label{E:tfi} \|\mT f\|_{H^{s_0}(\Gamma)} \leq s_j(i_1) \|\mT f\|_{H^{s}(\Gamma)}.\end{equation}
Combining (\ref{htf}), (\ref{E:hin}), and (\ref{E:tfi}), one arrives at  $[s_j(i_2)]^{1-\mu} \leq C_s [s_j(i_1)]^{\mu}$, where $C_s$ is a constant independent of $j$. Due to Lemma \ref{sing}, this implies $[j^{-s_1}]^{(1-\mu)} \leq C_s
 [j^{\frac{s_0-s}{n}}]^{\mu}$. Choosing $s$ big enough such that $s_1 (1-\mu) < \frac{(s-s_0)\mu}{n}$, we have a contradiction by letting $j\rightarrow \infty$. The proof is completed. \eproof

\section*{Appendix}

\begin{prop}\label{P:cont} Assume that $X,Y,Z$ are Banach spaces and $Z$ is  continuously imbedded into $Y$. Let $T: X \longrightarrow Y$ be a bounded linear operator such that $T(X) \subset Z$. Then $T: X \longrightarrow Z$ is also bounded.\end{prop}

{\bf Proof} Due to the closed graph theorem \cite{Ru73}, if it is sufficient to prove that $T$ has closed graph. That is, if
$\{x_k\}_k \subset X$ such that $x_k \rightarrow x$ in $X$ and $Tx_k \rightarrow z$ in $Z$, then $Tx=z$. In fact, since the imbedding $Z \rightarrow Y$ is continuous and $Tx_k \rightarrow z$ in $Z$, we have $Tx_k \rightarrow z$ in $Y$. On the other hand, since $x_k \rightarrow x$ in $X$ and $T:X \longrightarrow Y$ is continuous, we have $Tx_k \longrightarrow Tx$ in $Y$. Due to the uniqueness of limit of the sequence $\{Tx_k\}_k$ in $Y$, we have $Tx =z$. This finishes the proof. \eproof

\section*{Acknowledgments} This work was supported in part by NSF DMS grants 0604778 and 0715090, and the grant KUS-C1-016-04 from King Abdullah University of Science and Technology. The author expresses his gratitude to the NSF and KAUST for the support.

The author is grateful to Professors P. Kuchment, P. Stefanov, and G. Uhlmann for useful discussions and comments.

\end{document}